\documentclass{article}
\usepackage{latexsym,amssymb,amsmath}
\usepackage{amsfonts}

\newtheorem{theorem}{Theorem}[section]
\newtheorem{lemma}[theorem]{Lemma}
\newtheorem{proposition}[theorem]{Proposition}

\newtheorem{remark}[theorem]{Remark}
\newenvironment{proof}[1][Proof]{\textbf{#1.} }{\ \mathbb Rule{0.5em}{0.5em}}

\newcommand{\po}{{\mathbb P}}
\begin{document}

\title{Well posedness of a stochastic hyperviscosity-regularized 
3D Navier-Stokes equation}
\author{B. Ferrario\\     Dipartimento di Matematica - Universit\`a di Pavia}
\date{}
\maketitle

\section{Introduction}

We call  stochastic Navier--Stokes problem the following:
\begin{equation}\label{stochNS}
\left\{
\begin{array}{l}
 \dfrac{du}{dt} -\nu \Delta u+\left(  u\cdot\nabla\right)  u+\nabla p
         = f+ n
\\[2mm]
 \text{div }u   =0
\\
 u|_{t=0}   =u_{0}
\end{array}
\right.
\end{equation}
Here $u=u(t,x)$ is the 3-dimensional velocity vector field 
defined for $t\geq 0$ and
$x\in D\subseteq \mathbb R^3$, $p=p(t,x)$ is the
scalar pressure field, $\nu>0$ is the coefficient of kinematic
viscosity, 
$u_0=u_0(x)$ is the initial velocity,
$f=f(t,x)$ and $n=n(t,x)$ 
are, respectively, the
deterministic and stochastic forcing terms.
If the spatial domain has a boundary, we assume that $u$ vanishes on
$\partial D$.

When there is no noise term $n$, 
this reduces to the deterministic Navier--Stokes problem which models the motion
of viscous fluids.
For the 3-dimensional setting, both in the deterministic and in the stochastic case
we know the existence of a weak solution but uniqueness is proved in a smaller class, 
where existence is not known. 
The question of proving existence and uniqueness on any finite time interval 
and with any initial data for the deterministic
Navier--Stokes equation  is one of the  Millennium Prize problems (see, e.g., \cite{feff}).
However, for the 3-dimensional problem there are results for small initial data or locally in time,
whereas the 2-dimensional problem is well posed
(see, e.g., \cite{temam, temam2} for the deterministic problem and
\cite{Fla, funiq} for the 2-dimensional stochastic problem with
additive noise).

There have been many attempts to  modify 
the 3-dimensional Navier--Stokes equation
in order to prove an existence and uniqueness result. 
The first models go back to Lions \cite{lions}.
For more recent results, we focus on two particular
cases:  \cite{ms} and \cite{sri}.
In \cite{ms} 
the Laplacian operator $-\Delta$ is replaced with
$(-\Delta)^\alpha$ (for $\alpha >1$) in the {\it deterministic} Navier--Stokes equation;
the {\it stochastic} problem with a similar modification ($-\nu \Delta$ replaced with 
$-\nu_0\Delta+\nu_1(-1)^\alpha\Delta^\alpha$) is studied in \cite{sri}.
Setting $\alpha>1$ we obtain a model for
hyperviscous fluids (see \cite{sri} and references therein).

Our aim is to analyse the  well posedness of the
stochastic version of the modified Navier--Stokes equation
considered by Mattingly--Sinai in \cite{ms}, that is
\begin{equation}\label{alpha}
\left\{
\begin{array}{l}
 \dfrac{du}{dt}+ \nu (-\Delta)^\alpha u+\left(  u\cdot\nabla\right)  u+\nabla p
 =
  n
\\[2mm]
 \text{div }u   =0
\\
 u|_{t=0}   =u_{0}
\end{array}
\right.
\end{equation}
We shall consider the model of additive noise, i.e. $n$ is
independent of $u$. This is the simplest case, which reduces the
technicalities. However, the case of multiplicative noise can be
treated in a similar way.

In Section 3, we shall prove an existence and uniqueness
result for $\alpha \ge \frac 54$, as conjectured in \cite{f-gir}. 
The  bound  $\alpha > \frac 54$ appeared 
first in \cite{ms} for the deterministic problem. Regularity results
will be given in Section 4.

Finally, we point out that our technique is different 
from that of \cite{ms} or \cite{sri}; indeed, we use tools from
\cite{temam} and \cite{Fla}.

\section{Notations and preliminaries}        \label{prelim}
Let the spatial domain be a torus, i.e. 
the spatial variable $x$ belongs to 
$\mathcal T=[0,L]^3$ and
periodic boundary conditions are assumed.

We introduce the classical spaces for the Navier--Stokes equation
(see, e.g., \cite{temam2, temam} for all the results in this section).
$\mathcal{D}^{\infty}$ is defined 
as the space of infinitely differentiable
divergence free periodic fields $u:\mathcal{T} \to \mathbb R^3$, 
with zero mean ($ \int_{\mathcal{T}} u(x)  dx=0$).
Let $H^0$ be the closure of $\mathcal{D}^{\infty}$ 
in the $[L^2(\mathcal T)]^3$-topology; it is the subspace of $[L^2(\mathcal T)]^3$
of all fields $u$ such that $\mbox{div}\, u=0$,
the normal component of $u$  on the boundary is periodic, 
$\int_{\mathcal{T}}u\left(  x\right)  dx=0$.
We endow $H^0$
with the inner product
$
 \left\langle u,v\right\rangle=
\int_{\mathcal{T}} u(x)\cdot v(x)\ dx
$
and the associated norm $\left| \cdot \right|$.
\\
Similarly, for $m \in \mathbb N$ let
$H^m$ be the closure of $\mathcal{D}^{\infty}$ 
in the $[H^m(\mathcal T)]^3$-topology.

Let $A:D(A)\subset H^0\rightarrow H^0$ be the operator $Au=-\Delta u$
(componentwise) with $D(A)=H^2$.  This is called Stokes operator and it is
a strictly positive unbounded self-adjoint operator in $H^0$, whose
eigenvectors $h_j$ form a complete orthonormal basis
of the space $H^0$; the eigenvalues $\lambda_j$ are strictly positive
and
$0<\lambda_1\le \lambda_2 \le \dots$ with $\lambda_j \sim j^{2/3}$ for
$j \to \infty$.
Since the spatial domain is the torus, we know the
expressions of the eigenvectors with their eigenvalues
(see, e.g., \cite{FF}).

The power operators $A^\alpha$ are defined for any $\alpha \in \mathbb R$.
If $u=\sum_j u_j h_j$, then 
\[
A^\alpha u=\sum_{j=1}^\infty \lambda^\alpha_j u_j h_j
\quad\text{ and }\quad
|A^\alpha u|^2 = \sum_{j=1}^\infty \lambda_j^{2\alpha}|u_j|^2.
\]
We set $|u|_{2\alpha}=|A^\alpha u|$ and the space $H^s$ can be defined
(for any $s \in \mathbb R$) as the closure of $\mathcal{D}^{\infty}$ 
in the $|\cdot|_s$-metric.
For $s<0$ the space $H^s$ is the dual space of $H^{-s}$ with respect 
to the $H^0$-topology. 
The space
$H^{s+r}$ is dense and compactly embedded in $H^s$ for any $r>0$.

Notice that $|u|_m$ is equivalent to the usual 
$[H^m(\mathcal T)]^3$-norm.

In particular
\[
 |u|_1^2= \langle u, Au\rangle=\sum_{i,j=1}^{3}\int_{\mathcal T}
(\partial_j u_i(x))^2 dx.
  \]

We have $H^0 =\text{span} \{ h_{k}\} $
and we set $H_n=\text{span} \{h_{k}: |k|\le n\}$; 
moreover, we denote by
$\pi_n$ the projection operator from $H^0$ onto $H_n$.
The operators $A$ and $\pi_n$ commute.
By $\Pi$ we denote the projector operator from $[L^2(\mathcal T)]^3$ onto $H^0$.
\\
The operator $-A$ generates in $H^0$ 
(and in any $H^s$) an analytic semigroup
of negative type $\{e^{-tA}\}_{t\ge 0}$
of class $C_0$.

Let $B\left(  \cdot,\cdot\right) : H^1\times H^1\rightarrow H^{-1}$ be the
bilinear operator defined as
\begin{equation}\label{defB}
 \left\langle w,B\left(  u,v\right)  \right\rangle 
=
 \sum_{i,j=1}^{3}\int_{\mathcal T} u_i (\partial_i v_j) w_j dx
\end{equation}
for every $u,v,w\in H^1$. 
By the incompressibility condition, we have
\begin{equation}\label{incomp}
 \langle B\left(  u^1,u^2\right),u^2  \rangle =0, \qquad
 \langle B\left(  u^1,u^2\right),u^3  \rangle = 
       -\langle B\left(  u^1,u^3\right),u^2  \rangle .
\end{equation}
We shall use the following estimates (see Lemma 2.1 in \cite{temam2}):
\begin{equation}\label{Bcon4}
 |\langle B(u^1,u^2),u^3\rangle|
 \le \ c \  |u^1|\  |u^2|_\alpha  \  |u^3|_{\alpha} 
 \qquad \text{ for }\alpha\ge \frac 54,
\end{equation}
\begin{equation}\label{BconA}
 \begin{split}
 |\langle B(u^1,u^2),Au^3\rangle| 
 &\le 
 c |u^1|_\alpha |u^2|_1 |Au^3|_{\alpha-1}
  \qquad \text{ for }\alpha \ge \frac 54
 \\
 &=
 c |u^1|_\alpha |u^2|_1 |u^3|_{\alpha+1}
 \end{split}
\end{equation}
and similarly
\begin{equation}\label{BconA-2}
 \begin{split}
 |\langle B(u^1,u^2),Au^3\rangle| 
 &\le 
 c |u^1|_1 |u^2|_\alpha |Au^3|_{\alpha-1}
 \qquad \text{ for }\alpha \ge \frac 54
 \\
 &=
 c |u^1|_1 |u^2|_\alpha |u^3|_{\alpha+1}.
 \end{split}
\end{equation}
Here and in the following, we denote by $c$ a positive constant, which may vary 
from place to place.

\section{Main theorem}
We apply the projection operator $\Pi$ to the first equation in 
\eqref{alpha}. We get an It\^o equation in  an infinite 
dimensional Hilbert space:
\begin{equation} \label{sns} 
    \left\{
     \begin{array}{l}
     du(t) +\left[ \nu A^\alpha u(t)+B\big(u(t),u(t)\big)\:\right]\; dt
          =A^{-\gamma} \; dw(t)\\
     u(0)=u_0
     \end{array} 
    \right.
\end{equation}
assuming the noise is of white type in time and with spatial covariance independent of $u$.
This means that $w$ is a cylindrical Wiener process in $H^0$
defined on a complete probability space with filtration
$(\Omega,\mathcal F,\{\mathcal F_t\}_{t\ge 0},\mathbb P)$
(i.e. given a sequence $\{\beta_j\}_j$ of i.i.d. 
standard Wiener processes, we represent the Wiener process in series
as $w(t)= \sum_j \beta_j(t) h_j$).
For simplicity we consider the operator in front of $w$ to be a power of the Stokes operator;
this is the model studied in \cite{f-gir}.
\\
For $\alpha=1$, this corresponds to  the stochastic Navier--Stokes equation as analysed 
for instance in \cite{Fla} for the 2-dimensional setting. 
\\
The technique to study equation \eqref{sns} comes from 
\cite{bens, vishik, dpz, Fla}.
First  we consider the  linear
equation, that is
the modified stochastic Stokes equation
\begin{equation} \label{ou} 
    \left\{
     \begin{array}{l}
     dz(t) +\nu A^\alpha z(t) dt =A^{-\gamma} \; dw(t)\\
     z(0)=0
     \end{array} 
    \right.
\end{equation}
Then, we define $v:=u-z$; this unknown solves the following equation, obtained 
subtracting equation \eqref{ou} to equation \eqref{sns} 
and bearing in mind the bilinearity of the operator $B$: 
\begin{equation} \label{eq-v}
    \left\{
     \begin{array}{l}
     \dfrac{d\;}{dt}v(t) +\left[ \nu A^\alpha v(t)+B\big(v(t),v(t)\big) +
             B\big(z(t),v(t)\big) +B\big(v(t),z(t)\big) \:\right]
          \\\hspace*{3cm}=-B\big(z(t),z(t)\big) \\
     v(0)=u_0
     \end{array} 
    \right.
\end{equation}
The noise term $A^{-\gamma} \; dw(t)$ has disappeared.

Let $[0,T]$ be any finite time interval. We now state our main result.
\begin{theorem}\label{teo1}
Let $\alpha \ge \frac 54$.\\
For any $u_0\in H^1$, if $\gamma>\frac 34$ 
then there exists a unique process $u$ which is
a strong solution
of \eqref{sns} such that
$$
 u \in C([0,T];H^1)\cap L^{\frac {2\alpha}{\alpha-1}}(0,T;H^{\alpha}) \qquad \po-a.s.;
$$
$u$ is progressively measurable in these topologies and 
is a Markov process in $H^1$. 
\end{theorem}

\subsection{Existence}
We study pathwise the problems for the unknowns $z$ and $v$. This will imply an existence result 
for $u$.

For the linear problem we have (see, e.g., Proposition 4.1 in \cite{f-gir}, 
based on \cite{dpz})
\begin{lemma}\label{lemma-z}
If
\begin{equation} \label{conv-z}
 \alpha+2\gamma>\theta+\frac 32,
\end{equation}
then 
equation \eqref{ou} has a unique strong solution $z$ such that
\begin{equation}\label{nsz:2p}
 \po \{z\in C([0,T];H^{\theta})\}=1.
\end{equation}
\end{lemma}

Now, we work pathwise for the equation satisfied by $v$ and therefore also for $u$. 
\begin{proposition} \label{propo-v}
Let $\alpha \ge \frac 54$. 
\\
For any $u_0\in H^1$, if $\gamma>\frac 34$
then there exists a process $v$ which is a strong solution
of \eqref{eq-v} such that
$$
 v \in C([0,T];H^1)\cap L^2(0,T;H^{1+\alpha}) \qquad \po-a.s.
$$  
$v$ is progressively measurable in these topologies.
\end{proposition}
\proof
From Lemma \ref{lemma-z}
we have that  $z \in C([0,T];H^\alpha)$ $\po$-a.s., 
since $\gamma>\frac 34$.  
We take the scalar product of equation \eqref{eq-v} with $v$
and use \eqref{incomp}-\eqref{Bcon4}:
\[
\begin{split}
 \frac 12 \frac {d\;}{dt}|v|^2+\nu|v|^2_{\alpha}
 & = -  \langle B(v,v)+B(v,z)+B(z,v)+B(z,z), v\rangle 
\\
&=-
  \langle B(v,z)+B(z,z), v\rangle 
\\
 &\le c  |v| |z|_\alpha |v|_\alpha + c |z|_\alpha^2 |v|_\alpha
\\
 & \le \frac \nu 2 |v|_\alpha^2+c_\nu 
     \big(|z|^2_\alpha |v|^2 +|z|_\alpha^4\big) .
\end{split}
\]
Then
$$
 \frac {d\;}{dt}|v|^2
 \le
 c  |z|^2_\alpha |v|^2 +c |z|_\alpha^4
$$ 
and from Gronwall lemma: $\displaystyle\sup_{0\le t\le T} |v(t)|^2 <\infty$.
Moreover, integrating in time the first inequality above,
 we have $\int_0^T |v(t)|^2_\alpha dt <\infty$.

Now we take the scalar product of equation \eqref{eq-v} with $Av$:
$$\frac 12 \frac {d\;}{dt}|v|_1^2+\nu|v|^2_{1+\alpha}=-
  \langle B(v,v)+B(v,z)+B(z,v)+B(z,z), Av\rangle .
$$
We use \eqref{BconA} and Young inequality to obtain
\begin{equation} \label{st-per-v}
 \frac 12 \frac {d\;}{dt}|v|_1^2+\nu |v|^2_{1+\alpha}
 \le 
 \frac \nu 2 |v|^2_{1+\alpha} + c_\nu (|z|_\alpha^2+|v|_\alpha^2)|v|_1^2
  +c_\nu |z|_\alpha^4.
\end{equation}
As usual, from
$$
 \frac {d\;}{dt}|v|_1^2 \le 2 c_\nu (|z|_\alpha^2+|v|_\alpha^2)|v|_1^2
  + 2 c_\nu |z|_\alpha^4,
$$
Gronwall inequality, with the result $v \in L^2(0,T;H^\alpha)$ proved before, implies 
$\displaystyle \sup_{0\le t \le T}|v(t)|^2_1 <\infty$
and integrating in time \eqref{st-per-v}
we get $\int_0^T |v(t)|^2_{1+\alpha}dt <\infty$.

The technique to prove existence is classical (see \cite{temam}).
We consider first the 
finite dimensional problem in the unknown $v_n$, obtained 
projecting equation \eqref{eq-v} onto $H_n$. This is the Galerkin approximation, 
for any $n=1,2,\dots$.
The above estimates hold uniformly also for the Galerkin sequence: for
any $n$
$$
 \sup_{0\le t \le T}|v_n(t)|^2_1 <c_1, \qquad
  \int_0^T |v_n(t)|^2_{1+\alpha}dt < c_2
$$
for constants $c_1$ and $c_2$ independent of $n$.
\\
Any finite dimensional (Galerkin) problem has a solution.
By passing to the limit as $n \to \infty$ we get an existence result 
for \eqref{eq-v}.
We also need that 
$\frac{dv_n}{dt}$ is uniformly bounded in 
$ L^2(0,T;H^{1-\alpha})$.
We verify  easily that $\frac{dv_n}{dt}\in L^2(0,T;H^{1-\alpha})$ 
and the norm is bounded uniformly for all $n$,
since $A^\alpha v \in L^2(0,T;H^{1-\alpha})$ and according to \eqref{BconA} 
all the bilinear terms  in \eqref{eq-v}
are in $L^2(0,T;H^0)$. 
Therefore, we have a compact embedding (see Theorem 2.1, Ch. III
in \cite{temam}) so the Galerkin sequence stays in a compact subset of $L^2(0,T;H^1)$.
Therefore there exists a subsequence 
converging to $v$ as follows:
\[
\begin{split}
& v_m \to v  \quad \text{ weakly in }  L^2(0,T;H^{1+\alpha}),\\
& v_m \to v  \quad\star-\text{weakly in } L^\infty(0,T;H^1) ,\\
& v_m \to v   \quad\text{ strongly in } L^2(0,T;H^1).  
 \end{split}
\]
The strong convergence allows one 
to pass to the limit in the bilinear term (see Lemma 3.2, Ch. III
in \cite{temam}).
Finally, $v \in C([0,T];H^1)$ (see Lemma 1.2, Ch. III in \cite{temam}).
The limit $v$ fulfils all the estimates found above: 
$v\in C([0,T];H^1)\cap L^2(0,T;H^{1+\alpha})$.
$\hfill \Box$
\medskip

We conclude  noting that by interpolation
$L^\infty(0,T;H^1)\cap L^2(0,T;H^{1+\alpha})\subset
L^{\frac {2\alpha}{\alpha-1}}(0,T;H^\alpha)$.
Since the paths of $z$ belong to $C([0,T];H^\alpha)$ and 
those of $v$ to $L^\infty(0,T;H^1)\cap L^2(0,T;H^{1+\alpha})$,
then $u=v+z \in C([0,T];H^1)\cap L^{\frac {2\alpha}{\alpha-1}}(0,T;H^\alpha)$ 
$\po$-a.s. This concludes the existence result of Theorem \ref{teo1}
\\
The measurability property is inherited by from the Galerkin sequence.

\begin{remark}
The spatial covariance of the noise can be taken of a more general form. Indeed,
what is needed is that pathwise we have $z \in C([0,T];H^\alpha)$.
Hence, we can prove the same result of Theorem \ref{teo1} 
when instead of $A^{-\gamma} dw(t)$ the noise is $Gdw(t)$ 
assuming that the linear operator $G:H^0\to H^0$ is 
a Hilbert--Schmidt operator. 
This allows to consider any finite noise, that is acting on a finite number of components $h_k$
of the space $H^0$.
\\
More generally, Lemma \ref{lemma-z}
is true if the range of the operator $G$ is a subset of $D(A^\gamma)$
with $\gamma$ fulfilling \eqref{conv-z}.
\end{remark}

\subsection{Pathwise uniqueness} \label{unic}
We consider two solutions
$u_1$ and $u_2$ of \eqref{sns} obtained in the previous section;
we have that, for $\alpha \ge \frac 54$, 
$$u_1, u_2 \in C([0,T];H^1)
\cap L^2(0,T;H^\alpha) \; \po-a.s.
$$  
since $\frac {2\alpha}{\alpha-1}>2$. The difference 
$U=u_1-u_2$ satisfies
\begin{equation}
      \frac{d\;}{dt}U(t) +\nu A^\alpha U(t)+B\big(u_1(t),u_1(t)\big)
      -B\big(u_2(t),u_2(t)\big) = 0.
\end{equation}
Since the operator $B$ is bilinear, this becomes
\begin{equation}\label{diffU}
      \frac{d\;}{dt}U(t) +\nu A^\alpha U(t)+B\big(u_1(t),U(t)\big)
      +B\big(U(t),u_2(t)\big) = 0.
\end{equation}
Taking the scalar product of \eqref{diffU} with $AU$ in $H^0$, we get
$$
\frac 12 \frac{d\;}{dt}|U(t)|_1^2+\nu |U(t)|^2_{1+\alpha}=
 -\langle  B\big(u_1(t),U(t)\big)+B\big(U(t),u_2(t)\big), AU(t)\rangle 
$$ 
with $U(0)=0$.

We estimate the r.h.s. according to \eqref{BconA}-\eqref{BconA-2}
$$
 \frac 12 \frac{d\;}{dt}|U(t)|_1^2+\nu |U(t)|^2_{1+\alpha}
 \le c |u_1(t)|_\alpha |U(t)|_1  |U(t)|_{1+\alpha} + c
      |u_2(t)|_\alpha |U(t)|_1 |U(t)|_{1+\alpha} .
$$
Thus, by Young inequality:
$$
 \frac{d\;}{dt}|U(t)|_1^2
\le c \big(|u_1(t)|_\alpha^2+|u_2(t)|_\alpha^2 \big)|U(t)|_1^2
$$
and, by Gronwall inequality
\begin{equation}\label{sGro}
 |U(t)|_1^2\le |U(0)|_1^2 \; e^{\textstyle\int_0^t c
   (|u_1(s)|_\alpha^2+|u_2(s)|_\alpha^2 ) ds}.
\end{equation}
Then $U(t)=0$ for all $t \ge 0 $, because $U(0)=0$.
This means that pathwise we have $u_1(t)=u_2(t)$ for all $t \ge 0 $.

\begin{remark} 
i) The Markov property of $u$ comes from the same properties for the Galerkin approximations
and from the uniqueness result (see, e.g., \cite{Fla}). 
\\
ii) The pathwise estimate \eqref{sGro}
 implies also the Feller property
in $H^1$,
that is given a sequence of solutions $u^j$ with initial data $u^j_0$, 
if $\displaystyle\lim_j u_0^j=u_0$ in $H^1$ 
then $\displaystyle\lim_j \mathbb E \phi (u^j(t))=\mathbb E \phi(u(t))$
for any $t \in [0,T]$ and any
 continuous bounded function $\phi:H^1 \to\mathbb R$.
\end{remark}

\section{Regularity results}
Considering as phase space other Hilbert spaces $H^s$, 
we get different bounds on $\alpha$ to obtain that the dynamics 
of the stochastic Navier--Stokes equation is well posed in such $H^s$.
As it has been pointed out in \cite{ms}, the smaller is $s$ 
(with $u_0\in H^s$) the bigger is $\alpha$.
Theorem \ref{teo1}  deals with  $s=1$.
In this section, we show  that problem \eqref{sns} is well-posed in
the space $H^0$ if $\alpha>\frac 32$, and
in the spaces $H^s$ with $s\ge 2$ if $\alpha \ge
  \frac 54$.
\\[2mm]
$\boxed{\mathbf H^0\text{\bf -regularity}}$\\
We have the following result
\begin{proposition}
Let $\alpha>\frac 32$.
\\
For any $u_0\in H^0$, if $\gamma>\frac 34$
then there exists 
a unique process $u$ which is
a strong solution
of \eqref{sns} such that
$$
 u \in C([0,T];H^0)\cap 
L^2(0,T;H^{\alpha}) \qquad \po-a.s.;
$$
$u$ is progressively measurable in these topologies 
and a Markov process in $H^0$. 
\end{proposition}
\proof
Existence is proved by means of a priori estimates as in the proof 
of Proposition
\ref{propo-v}; to be precise, for $\alpha \ge \frac 54$ we get that
there exists a solution $u \in C([0,T];H^0)\cap
L^2(0,T;H^\alpha)$ $\po$-a.s.
if $z$ has paths in $C([0,T];H^\alpha)$, i.e.
if $\gamma > \frac 34$.

Pathwise 
uniqueness is obtained according to the result by Prodi
\cite{prodi}, requiring $u \in L^{s}(0,T;[L^q(\mathcal T)]^3) \ \po$-a.s.
for $\frac 2s +\frac 3 q \le 1$.
However, using an interpolation result and Sobolev embedding we have
$$
 L^\infty(0,T;H^0)\cap L^2(0,T;H^\alpha) \subset L^s(0,T;H^{2\frac \alpha s})
 \subset L^s(0,T;[L^q(\mathcal T)]^3) 
$$
for $2< s < \infty$ and $\frac 1q=\frac 12 -\frac {2\alpha}{3s}$. 
The condition $\frac 2s +\frac 3 q \le 1$ 
holds if $\alpha>\frac 32$. $\hfill \Box$

\begin{remark}
For $\alpha \ge 1$ we can prove 
 that there exists a solution of equation \eqref{sns} such that
$u \in C_w([0,T];H^0)\cap L^\infty(0,T;H^0)\cap
L^2(0,T;H^\alpha)$ $\po$-a.s. But uniqueness is unknown.

Consider, for instance, $\alpha=1$.  
We require that  
$z \in C([0,T];H^{\frac 32})$ $\po$-a.s. and equation \eqref{eq-v} 
is treated as in the deterministic setting.\\
For this, change the first estimate in the proof of Proposition 
\ref{propo-v} as follows:
\[
\begin{split}
\frac 12 \frac{d}{dt}|v|^2+\nu|v|_1^2
 & =-\langle B(v+z,z),v\rangle
\\
 &\le c |v+z| |z|_{\frac 32} |v|_1 
\\
 & \le c |v| |z|_{\frac 32} |v|_1+ |z|^2_{\frac 32} |v|_1 
\\
 &\le \frac \nu 2 |v|_1^2
 + c_\nu  |z|_{\frac 32}^2 |v|^2 +  c_\nu  |z|_{\frac 32}^4.
\end{split}
\]
Then $\sup_{0\le t \le T} |v(t)|<\infty$, 
$\int_0^T |v(t)|_1^2 dt<\infty$.

Moreover, from Lemma 2.1 in \cite{temam2} we have that
 $B:H^0\times H^1 \to H^{-\sigma}$ for any $\sigma
>\frac 32$. Then
$\dot v= -\nu A v-B\big(v+z,v+z\big)
\in L^2(0,T;H^{-\sigma})$. This gives a compact embedding and 
therefore there exists a subsequence of the Galerkin sequence
that converges to $v$ as follows:
\[
\begin{split}
& v_m \to v  \;\;\text{ weakly in }  L^2(0,T;H^1),\\
& v_m \to v \;\;\star-\text{weakly in } L^\infty(0,T;H^0) ,\\
& v_m \to v  \;\;\text{ strongly in } L^2(0,T;H^0).  
 \end{split}
\]
Finally $v \in C_w([0,T];H^0)$.
Notice that the previous result $v \in C([0,T];H^0)$ came from
$v \in L^2(0,T;H^\alpha), \dot v \in L^2(0,T;H^{-\alpha})$
(see Lemma 1.2, Ch. III in \cite{temam}).
\end{remark}
\begin{remark}
To compare our result with \cite{sri}, we have that
\cite{sri}, for its model, deals with the phase space $H^0$
assuming $\alpha\ge 2$.
\end{remark}

\pagebreak
\noindent
$\boxed{\mathbf {H^s\text{\bf -regularity with } s \ge 2}}$\\
We need the following estimates:
\begin{lemma}
\begin{equation}\label{B1}
|B(u,\tilde u)|_m \le c |u|_{m+1}|\tilde u|_{m+1} 
\qquad \text{ for } m =1,2,3,\dots
\end{equation}
\end{lemma}
\proof
First, consider \eqref{B1} for $m=1$.
We have 
\[
\begin{split}
 |B(u,\tilde u)|_1^2 &=|(u\cdot \nabla)\tilde u|_1^2 
\\&=
\sum_{k,l=1}^3\big|\partial_k(\sum_{i=1}^3 u_i\partial_i\tilde u_l)\big|_{L^2}^2 
\\
 &\le 2 \sum_{k,l=1}^3\big|\sum_{i=1}^3 \partial_k  u_i\partial_i \tilde u_l\big|_{L^2}^2
 + 2 \sum_{k,l=1}^3\big|\sum_{i=1}^3 u_i\partial_k\partial_i \tilde u_l\big|_{L^2}^2
\\
&\le 6 \sum_{k,l,i}\big| \partial_k  u_i\big|^2_{L^4} 
                 \big|\partial_i \tilde u_l\big|_{L^4}^2
 + 6 \sum_{k,l,i}\big| u_i\big|_{L^\infty}^2  
                \big|\partial_k\partial_i \tilde u_l\big|_{L^2}^2.
\end{split}\]
Then use the continuous embeddings 
$H^1(\mathcal T) \subset L^4(\mathcal T)$ and
$H^2(\mathcal T)\subset L^\infty(\mathcal T)$.

For $m=2,3,\dots$ we use that 
$H^m$ is a multiplicative algebra if $m >\frac 32$; then
\[
|B(u,\tilde u)|_m \le c |u|_m|\tilde u|_{m+1} \text{ for } m =2,3,\dots
\]
which is even stronger than \eqref{B1}.
$\hfill \Box$

\medskip
We have the following result
\begin{proposition}
Let $\alpha\ge \frac 54$ and $s\ge 2$.
\\
For any $u_0\in H^s$, if $\alpha+2\gamma>s+\frac 32$
then there exists 
a unique process $u$ which is
a strong solution
of \eqref{sns} such that
$$
u \in C([0,T];H^s) 
$$
$\po$-a.s.\\
$u$ is progressively measurable in these topologies 
and is a Markov process in $H^s$. 
\end{proposition}

For simplicity,
we provide the proof for  $s=2$. In this way we show the
difference with respect to the case $s=1$ considered in the previous
section. However, the proof would go along the same lines for
$s>2$ using \eqref{B1}.

\proof Set $s=2$. Then 
almost every  path of $z$ is in $C([0,T];H^2)$.
\\
We prove existence for $\alpha \ge \frac 54$.
We use \eqref{B1}  with $m=1$ and
take the scalar product of equation \eqref{eq-v} with $A^2v$:
\[\begin{split}
\frac 12 \frac {d\;}{dt}|v|_2^2+\nu|v|^2_{2+\alpha}
=-
  \langle B(v+z,v+z), A^2v\rangle 
&
= -\langle A^{\frac 12}B(v+z,v+z), A^{\frac32}v\rangle \\
&
\le |B(v+z,v+z)|_1 |v|_3\\
&
\le c |v+z|_2^2 |v|_3\\
&
\le c |v+z|_2^2 |v|_{2+\alpha}\\
&\le
 \frac \nu 2 |v|^2_{2+\alpha}+ c_\nu |v|_2^4+c_\nu |z|^4_2.
\end{split}\]
Since we already know from Proposition \ref{propo-v}
that $v \in L^2(0,T;H^{1+\alpha})
\subset L^2(0,T;H^2)$, it follows as usual by Gronwall lemma that
$v \in L^\infty(0,T;H^2)\cap L^2(0,T;H^{2+\alpha})$. From now on, the proof
goes as in Proposition \ref{propo-v}.
\\
Pathwise uniqueness: the estimates hold for any $\alpha \ge 1$ but the
regularity required on $u_i$ holds for $\alpha \ge \frac 54$. This
shows that in $H^2$ it is ''easier'' to prove uniqueness  than existence.
\\
Set $U=u_1-u_2$ as in Section \ref{unic}; now $u_1, u_2 \in C([0,T];H^2)$.
Taking the scalar product of \eqref{diffU} with $A^2U$ in $H^0$, we get
$$
\frac 12 \frac{d\;}{dt}|U(t)|_2^2+\nu |U(t)|^2_{2+\alpha}=
 -\langle  B\big(u_1(t),U(t)\big)+B\big(U(t),u_2(t)\big), A^2U(t)\rangle 
$$ 
with $U(0)=0$.
As before, we estimate the r.h.s. by means of \eqref{B1}, and get
\[\begin{split}
 \frac 12 \frac{d\;}{dt}|U(t)|_2^2+\nu |U(t)|^2_{2+\alpha}
 &\le c |u_1(t)|_2 |U(t)|_2  |U(t)|_3 + c
      |u_2(t)|_2 |U(t)|_2 |U(t)|_3 
\\
&\le
 c |u_1(t)|_2 |U(t)|_2  |U(t)|_{2+\alpha} + c
      |u_2(t)|_2 |U(t)|_2 |U(t)|_{2+\alpha}
\\
&\le
 \frac \nu 2 |U(t)|_{2+\alpha}^2+
 c_\nu(|u_1(t)|_2^2+|u_2(t)|_2^2)|U(t)|_2^2. 
\end{split}
\]
From 
$$
\frac{d\;}{dt}|U(t)|_2^2\le  2 c_\nu(|u_1(t)|_2^2+|u_2(t)|_2^2)|U(t)|_2^2
$$
we conclude that $|U(t)|_2=0$ for any $t \in [0,T]. \hfill \Box$

\end{document}